\documentclass[11pt,reqno]{amsart}
\usepackage{amssymb}
\usepackage{amsfonts}
\usepackage{amsmath}
\usepackage{graphicx}
\usepackage{hyperref}
\usepackage{float}
\usepackage{epstopdf}
\usepackage{color}
\usepackage{bm}
\usepackage{comment}
\usepackage{soul}

\allowdisplaybreaks
\newtheorem{theorem}{Theorem}[section]
\newtheorem{proposition}[theorem]{Proposition}
\newtheorem{lemma}[theorem]{Lemma}

\newtheorem{definition}[theorem]{Definition}

\def\SG{\mathcal{SG}}

\def\lSG{\widetilde{\SG}}

\numberwithin{equation}{section}

\begin{document}
\title[Spectral analysis beyond $\ell^2$ on Sierpinski lattices]{Spectral analysis beyond $\ell^2$ on Sierpinski lattices}

\author{Shiping Cao}
\address{Department of Mathematics, Cornell University, Ithaca 14853, USA}
\email{sc2873@cornell.edu}
\thanks{}

\author{Yiqi Huang}
\address{Department of Mathematics,	The Chinese University of Hong Kong, Shatin, N.T., Hong Kong}
\email{1155091910@link.cuhk.edu.hk}
\thanks{}

\author{Hua Qiu}
\address{Department of Mathematics, Nanjing University, Nanjing 210093, China}
\email{huaqiu@nju.edu.cn}
\thanks{The research of Qiu was supported by the NSFC grant 12071213}

\author{Robert S. Strichartz}
\address{Department of Mathematics, Cornell University, Ithaca 14853, USA}
\email{str@cornell.math.edu}

\author{Xiaohan Zhu}
\address{Department of Mathematics, The University of Wisconsin-Madison, Madison 53706, USA}
\email{xzhu274@wisc.edu}
\thanks{}

\subjclass[2010]{Primary 47A10, 05C63, 28A80}

\date{}

\keywords{Sierpinski lattices, Laplacian, spectrum, eigenfunctions}

\begin{abstract}
We study the spectrum of the Laplacian on the Sierpinski lattices. First, we show that the spectrum of the Laplacian, as a subset of $\mathbb{C}$, remains the same for any $\ell^p$ spaces. Second, we characterize all the spectral points for the lattices with a boundary point.  
 \end{abstract}

\maketitle
\section{introduction}
In this note, we study the spectrum of the Laplacian on the Sierpinski Lattice $\lSG$. The problem was fully investigated by A. Teplyaev \cite{T} in the $\ell^2$ setting. We continue his study for the $\ell^p$ case. The $l^p$-spectra of the Laplacian on graphs have been studied in a variety of works \cite{BHK,Davies2,Davies1} in history. In particular, in \cite{BHK}, for graphs of uniform sub-exponential volume growth, which include the Sierpinski type lattices as particular examples, the $p$-independence of the spectra of the Laplacian has been shown. In this work, under the more concrete setting, we will try to look more carefully into the spectral points. In particular, the $l^1$-spectrum is of special interest for its different behavior.

\begin{figure}[htp]
	\includegraphics[width=8cm]{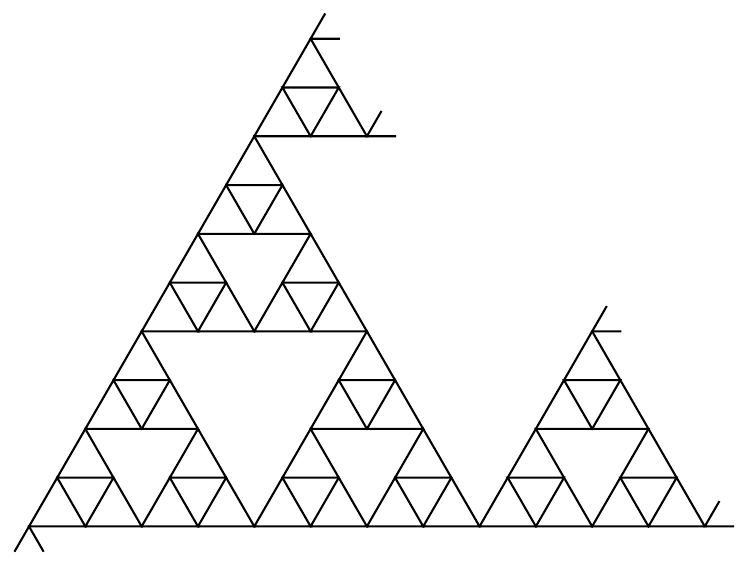}
	\caption{The Sierpinski lattices.}\label{fig1}
\end{figure}

The \textit{Sierpinski lattice} is an infinite graph defined as follows. Let 
$$F_0(x)=\frac{1}{2}x+(\frac{1}{4},\frac{\sqrt{3}}{4}),\quad F_1(x)=\frac{1}{2}x,\quad F_2(x)=\frac{1}{2}x+(\frac{1}{2},0),$$
and $q_0=(\frac{1}{2},\frac{\sqrt{3}}{2}),q_1=(0,0),q_2=(1,0)$ be the fixed points of $F_i$ respectively. Given an infinite word $\omega=\omega_1\omega_2\cdots \in \{0,1,2\}^\infty$, the corresponding Sierpinski lattice $\lSG_\omega$ (we will write $\lSG_\omega=\lSG$ when this causes no confusion) is constructed as
\[\lSG=\bigcup_{m=1}^\infty F_{\omega_1}^{-1}F_{\omega_2}^{-1}\cdots F_{\omega_m}^{-1}(V_m),\]
where $V_0=\{q_0,q_1,q_2\}$ and $V_m=\bigcup_{i=0}^2 F_i(V_{m-1})$ are defined iteratively. See Figure \ref{fig1}.
For $x,y\in \widetilde{\SG}$, we write $x\sim y$ if $x,y\in F_{\omega_1}^{-1}F_{\omega_2}^{-1}\cdots F_{\omega_m}^{-1}F_{l_m}F_{l_{m-1}}\cdots F_{l_1}(V_0)$ for a sequence $l_1,\cdots, l_m\in\{0,1,2\}$,
  and call $y$ a \textit{neighbouring vertex} of $x$. All vertices in  $\lSG$ have four neighbouring vertices, except  the \textit{boundary vertex}, which admits only two neighbouring vertices. The structure of $\lSG$ depends on the infinite word $\omega$, although we omit this dependence in our notation. In fact,  depending on the choice of $\omega$, there are uncountably many distinct Sierpinski lattices without boundary, but essentially only one with a non-empty boundary. The boundary exists (consists of one point) if and only if there exists $M\in \mathbb{N}$ and $i\in \{0,1,2\}$ such that $\omega_m=i$ for all $m\geq M$, see \cite{T}.

The Laplacian $\Delta$ (imposing the Neumann boundary condition at the boundary) on $\lSG$ is defined as
\begin{equation}\label{eqn11}
\Delta f(x)=\begin{cases}
\sum_{y\sim x}f(y)-4f(x), \quad &\text{if $x$ is not a boundary point,}\\ 
\sum_{y\sim x}2f(y)-4f(x),\quad &\text{if $x$ is a boundary point.}
\end{cases}
\end{equation}

The celebrated result of A. Teplyaev \cite{T} showed that the Laplacian $\Delta$, viewed as an operator $\ell^2(\lSG)\to\ell^2(\lSG)$ has pure point spectrum. What's more, $\ell^2(\lSG)$ admits a basis of localized eigenfunctions of the Laplacian, which can be generated using the spectral decimation method \cite{FS,tS}. 

Roughly speaking, the spectral decimation algorithm reveals a simple relationship between eigenvalues and eigenfunctions on approximating graphs (or fractal lattices) of different levels. In addition, the Neumann and Dirichlet eigenvalues naturally appears to be descendants of forbidden eigenvalues, which are special eigenvalues on the approximating graphs that the decimation formula fails to hold. Readers can find an easy introduction to the spectral decimation method in the book \cite{S1} (Chapter 3).

In particular, all the localized eigenfunctions and eigenvalues on $\lSG$ are generated with the spectral decimation method. What's more, as an application of the method, the relation between the spectral analysis on an infinite Sierpinski gasket, which is a particular example of fractal blowups introduced by R.S. Strichartz \cite{S0}, and on its corresponding lattice can be clearly described. See \cite{T} for a detailed discussion. For the spectral analysis on other fractal graphs and related fractalfolds, see \cite{BCDEHKMST1, BCDEHKMST2,HSTZ,MT,S2,ST}. We list some important references on analysis on fractals which base on or relate to the spectral analysis, \cite{BS,DRS,FKS,IRS,ORS,RS,S-1,S3,S4}.

In this work, we will describe the spectrum of the Laplacian $\Delta:\ell^p(\lSG)\to\ell^p(\lSG)$ for $1\leq p\leq\infty$. In Section 2, we will show that the $\ell^p$-spectrum is a union of a Julia set and a discrete set named $6$-series eigenvalues. The spectrum remains the same for any $1\leq p\leq \infty$ and any lattice by an application of spectral decimation method. On the other hand, in Section 3, we will see that for $1<p<\infty$, the Laplacian has only point spectrum and continuous spectrum, while $\Delta:\ell^1\to\ell^1$ has all three kinds of spectral points. This phenomenon is a consequence of the existence of the $4$-eigenfunctions in $\ell^\infty(\lSG)$, which does not live in other $\ell^p(\lSG)$ spaces. In addition, we get a full description of the $\ell^1$-spectrum for the lattices with one boundary vertex.

\section{The $\ell^p$-spectrum of $\Delta$ on $\lSG$}\label{intro}
In this section, we compute the spectrum  of the Laplacian 
$\Delta:\ell^p(\lSG)\to \ell^p(\lSG),$ 
which is stated in the following Theorem \ref{Q1}. We first remark that, as shown in \cite{T}, the Julia set $\mathcal{J}$ corresponding to $R(\lambda)=\lambda(5-\lambda)$ is contained in the spectrum of the $l^2$-Laplacian. For convenience of readers, we provide the definition of Julia sets as follows. Detailed discussions can be found in \cite{M}.

\begin{definition}
Let $f:\hat{\mathbb{C}}\to \hat{\mathbb{C}}$ be a non-constant holomorphic mapping, where $\hat{\mathbb{C}}=\mathbb{C}\cup\{\infty\}$ is the extended complex plane. $z\in \hat{\mathbb{C}}$ is called a \emph{normal point} if there is a neighbour $U$ of $z$ such that $\{f^{\circ{k}}\}_{k\geq 1}$ is a normal family, which means every infinite sequence of maps from $\{f^{\circ k}\}_{k\geq 1}$ contains a subsequence which
converges locally uniformly to a limit on $U$.  

The \emph{Fatou set} of $f$ consists of all normal points of $f$, and the \emph{Julia set} of $f$ is the complement of the Fatou set.
\end{definition}

Since $R$ is a quadratic polynomial, a well-known equivalent definition for $\mathcal{J}$ can be given as
\[\mathcal{J}=\{x\in\mathbb{C}: \{R^{\circ k}(x)\}_{k=0}^\infty\in l^\infty\}.\]
In other words, the Julia set $\mathcal J$ of $R$ consists of points in $\mathbb{C}$ with bounded orbit.

\noindent\textbf{Remark.} The Hausdorff dimension of $\mathcal{J}$ is the unique zero of the Bowen's function 
$$B(t)=\lim_{k\to\infty} \frac{1}{k}\ln \sum_{w\in R^{-k}(x)}|(R^{\circ k})'(w)|^{-t},$$
where $x$ is any chosen point in $\mathcal{J}$. See Section 9.1 in the book \cite{PU}.

\begin{theorem}\label{Q1}
	Let $\sigma(\Delta)=\mathcal{J}\cup \Sigma_6$, where $\Sigma_6=\{6\}\cup (\bigcup_{m=0}^\infty R^{\circ -m}\{3\})$. The spectrum of $\Delta: \ell^p(\widetilde{\SG})\rightarrow \ell ^p (\widetilde{\SG})$ is equal to $\sigma(\Delta)$ for all $1\leq p\leq \infty$.
\end{theorem}

\noindent\textbf{Remark.} Theorem 2.2 is proved for the $\ell^2(\lSG)$ case in \cite{T}. Here a different approach will be used to deal with general $\ell^p$ cases. The theorem is also valid for $\Delta:C_0(\lSG)\to C_0(\lSG)$.\vspace{0.2cm} 

We will prove Theorem \ref{Q1} with several lemmas.  For a fixed infinite word $\omega$, we consider a sequence of sparse lattices $\lSG^{(-k)}$ defined as
\[\lSG^{(-k)}=\bigcup_{m=k}^\infty F_{\omega_1}^{-1}F_{\omega_2}^{-1}\cdots F_{\omega_m}^{-1}V_{m-k},\]
and we say $x\sim_{-k}y$ if  $x,y\in F_{\omega_1}^{-1}F_{\omega_2}^{-1}\cdots F_{\omega_m}^{-1}F_{l_{m-k}}F_{l_{m-k-1}}\cdots F_{l_1}(V_0)$ for a sequence $l_1,\cdots, l_{m-k}\in\{0,1,2\}$. The Laplacian $\Delta_{(-k)}$ on $\lSG^{(-k)}$ can be defined in a similar manner as (\ref{eqn11}).
The eigenvalues and eigenfunctions on $\lSG^{(-k)}$ for different $k$'s are related by the spectral decimation method \cite{FS,ST,T}.

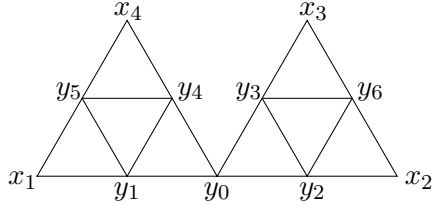
\begin{figure}[htp]
	\includegraphics[width=5cm]{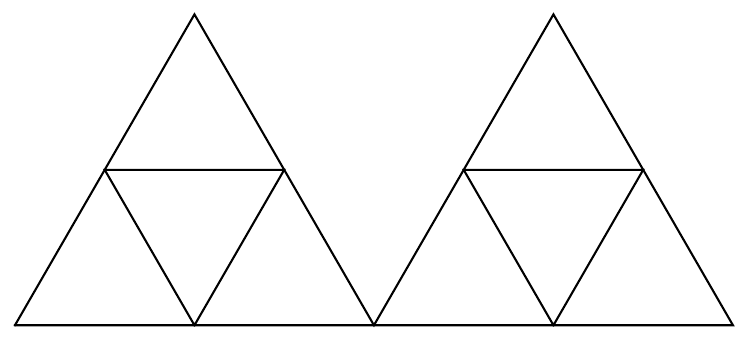}
	\begin{picture}(0,0)
	\put(-4,0){$x_2$}
	\put(-154,0){$x_1$}
	\put(-114,64){$x_4$}
	\put(-44,64){$x_3$}
	
	\put(-80,-4){$y_0$}
	\put(-68,34){$y_3$}
	\put(-44,-4){$y_2$}
	\put(-114,-4){$y_1$}
	\put(-90,34){$y_4$}
	\put(-136,34){$y_5$}
	\put(-22,34){$y_6$}
	\end{picture}
	\caption{The graph $\Gamma$.}\label{fig2}
\end{figure}

To understand the spectral decimation, we only need to focus on a small neighbourhood of a point $y_0$ in $\lSG^{(-k-1)}$. For convenience, we only consider the case $k=0$ and a point $y_0\in \lSG^{(-1)}$ with four neighbouring vertices, noticing that the boundary vertices and $k\geq 1$ cases can be dealt with in an essentially same way. Let $\{x_i\}_{i=1}^4$ be the neighbouring vertices of $y_0$ in $\lSG^{(-1)}$, and let $\{y_i\}_{i=1}^6$  be the vertices in $\widetilde{\SG}$ bounded by $\{x_i\}_{i=1}^4$. The induced subgraph in $\widetilde{\SG}$ is denoted by $\Gamma$ in the following context, see Figure \ref{fig2}. Clearly, the definition of $\Delta$ and $\Delta_{(-1)}$ on $\Gamma$ are naturally inherited from those on the graphs $\lSG$ and $\lSG^{(-1)}$ as follows,
\begin{eqnarray}
\label{eqn21}\Delta f(y_i)=\sum_{z\sim y_i}f(z)-4f(y_i),&0\leq i\leq 6,\\
\label{eqn22}\Delta_{(-1)} f(y_0)=\sum_{i=1}^4f(x_i)-4f(y_0).&
\end{eqnarray}

\begin{proposition}[Spectral decimation]\label{decimation}
Let $\lambda\notin \{2,5,6\}$ and $R(\lambda)=\lambda(5-\lambda)$. Let $\Delta$ and $\Delta_{(-1)}$ on $\Gamma$ be defined in (\ref{eqn21}) and (\ref{eqn22}). 

(a). Let $f\in l(\Gamma)$ and $-\Delta f(y_i)=\lambda f(y_i),\forall 0\leq i\leq 6$. Then $-\Delta_{(-1)}f(y_0)=R(\lambda)f(y_0)$. 

(b). Given any values $f(x_i),i=1,2,3,4$ and $f(y_0)$ such that $-\Delta_{(-1)}f(y_0)=R(\lambda)f(y_0)$, there is a unique extension $f\in l(\Gamma)$ such that $f$ satisfies the eigenvalue equations $-\Delta f(y_i)=\lambda f(y_i),\forall 0\leq i\leq 6$.
\end{proposition}
\noindent\textit{Proof.} This can be done by elementary calculation: first evaluate equation (\ref{eqn21}) at $y_1,y_4,y_5$ to yield $y_1,y_4$ as functions of $y_0,x_1,x_4,\lambda$, and do similarly to yield $y_2,y_3$ as functions of $y_0,x_2,x_3,\lambda$; then evaluate equations (\ref{eqn21}) and (\ref{eqn22}) separately at $y_0$ and compare the outcome values of $y_0$ in terms of $x_1,x_2,x_3,x_4,\lambda$; the \textit {forbidden eigenvalues} $2,5,6$ will emerge in the calculation. See Section 3.2 in the book \cite{S1} for details.\hfill$\square$\vspace{0.2cm}

Proposition \ref{decimation} can be easily applied to the eigenvalue problems on the Sierpinski lattices. However, to deal with the spectrum, we need somewhat stronger versions. We will do this in the following two lemmas.

\begin{figure}[htp]
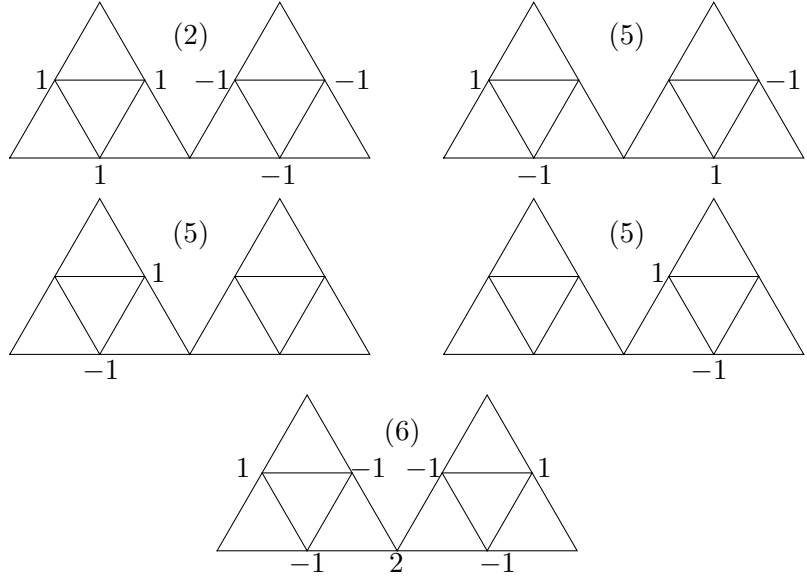

	\includegraphics[width=5cm]{gammat.pdf}\qquad
	\includegraphics[width=5cm]{gammat.pdf}\vspace{0.3cm}
	\includegraphics[width=5cm]{gammat.pdf}\qquad
	\includegraphics[width=5cm]{gammat.pdf}\vspace{0.3cm}
	\includegraphics[width=5cm]{gammat.pdf}
	\begin{picture}(0,0)
	\put(-212,178){$1$}\put(-167,178){$1$}\put(-190,142){$1$}\put(-160,195){(2)}
	\put(-152,178){$-1$}\put(-99,178){$-1$}\put(-127,142){$-1$}\put(5,195){(5)}
	\put(-48,178){$1$}\put(64,178){$-1$}\put(-30,142){$-1$}\put(43,142){$1$}\put(5,120){(5)}\put(-160,120){(5)}
	\put(-168,105){$1$}\put(-194,68){$-1$}
	\put(21,105){$1$}\put(36,68){$-1$}
	\put(-78,-5){$2$}\put(-80,45){(6)}
	\put(-72,31){$-1$}
	\put(-44,-5){$-1$}
	\put(-116,-5){$-1$}
	\put(-93,31){$-1$}
	\put(-136,31){$1$}
	\put(-22,31){$1$}

	\end{picture}
	\caption{The $2,5,6$-Dirichlet eigenfunctions on $\Gamma$ (with only non-zero values marked).}\label{fig3}
\end{figure}

\begin{figure}[htp]
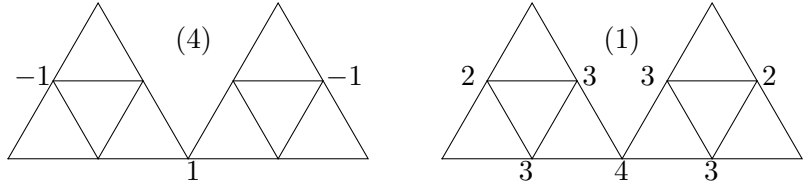

	\includegraphics[width=5cm]{gammat.pdf}\qquad
	\includegraphics[width=5cm]{gammat.pdf}
	\begin{picture}(0,0)
	\put(-78,-5){$4$}\put(-82,45){(1)}\put(-244,45){(4)}
	\put(-68,31){$3$}
	\put(-44,-5){$3$}
	\put(-114,-5){$3$}
	\put(-90,31){$3$}
	\put(-136,31){$2$}
	\put(-22,31){$2$}
	
	\put(-240,-5){$1$}
	\put(-305,31){$-1$}
	\put(-187,31){$-1$}
	\end{picture}
	\caption{The $4,1$-Dirichlet eigenfunctions on $\Gamma$.}\label{fig4}
\end{figure}

\begin{lemma}\label{Q1lemma1}
	Consider the Dirichlet eigenvalue problem on $\Gamma$:
	\[\begin{cases}
	f(x_i)=0,&1\leq i\leq 4,\\
	-\Delta f(y_i)=\lambda f(y_i),&1\leq i\leq 6.
	\end{cases}\]
	 All the Dirichlet eigenvalues are $\{1,2,4,5,6\}$.
\end{lemma}
\noindent\textit{Proof.} We can easily find one $2$-eigenfunction, three $5$-eigenfunctions and one $6$-eigenfunction, as shown in Figure \ref{fig3}, see Section 3.2, 3.3 in \cite{S1}. In addition, we find one $4$-eigenfunction and one $1$-eigenfunction as shown in Figure \ref{fig4}. All the above  give $7$ linearly independent eigenfunctions. \hfill$\square$

\begin{lemma}\label{Q1lemma2} Let $\lambda\notin\{1,2,4,5,6\}$. 
There exist constants $\{c_{i,\lambda}\}_{i=0}^6$ such that for any $f\in \l(\Gamma)$ we have 
$$(R(\lambda)+\Delta_{(-1)})f(y_0)=\sum_{i=0}^6 c_{i,\lambda}(\lambda+\Delta)f(y_i).$$
In addition, $c_{0,\lambda}\neq 0$.
\end{lemma}

\noindent\textit{Proof.} Let $f\in l(\Gamma)$. By the assumption and using Lemma \ref{Q1lemma1}, we know that $\lambda$ is not a Dirichlet eigenvalue. So there is a unique solution to each of the following boundary value problems.
\begin{equation}\label{Q1eqn1}
\begin{cases}
u(x_i)=0,&1\leq i\leq 4,\\
(\lambda+\Delta) u(y_i)=(\lambda+\Delta) f(y_i),&0\leq i\leq 6,
\end{cases}
\end{equation}
\begin{equation}
 \begin{cases}
 v(x_i)=f(x_i),&1\leq i\leq 4\\
 (\lambda+\Delta) v(y_i)=0,&0\leq i\leq 6.
 \end{cases}
\end{equation}
Clearly, we have
\[f=u+v.\]
In addition, since $v$ is an eigenfunction of $\Delta$, we have 
$$(R(\lambda)+\Delta_{(-1)})v(y_0)=0,$$
by using Proposition \ref{decimation}. As a consequence, we get $(R(\lambda)+\Delta_{(-1)})f(y_0)=(R(\lambda)+\Delta_{(-1)})u(y_0)$. On the other hand, since $u$ is uniquely determined by the linear equations (\ref{Q1eqn1}), we conclude that there are constants $c_{i,\lambda}$ such that 
\[(R(\lambda)+\Delta_{(-1)})f(y_0)=(R(\lambda)+\Delta_{(-1)})u(y_0)=\sum_{i=0}^6c_{i,\lambda}(\lambda+\Delta) f(y_i).\]

Lastly, we prove $c_{0,\lambda}\neq 0$ by contradiction. Assume $c_{0,\lambda}=0$, then we have 
\[f(y_0)=\big(R(\lambda)-4\big)^{-1}\big(-\sum_{i=1}^4f(x_i)+\sum_{i=1}^6 c_{i,\lambda}(\lambda+\Delta)f(y_i)\big).\]
In addition, it is easy to see that $\{f(y_i)\}_{i=1}^6$ are uniquely determined by  $\{f(x_i)\}_{i=1}^4,\{(\Delta+\lambda)f(y_i)\}_{i=1}^6$ and $f(y_0)$. As a consequence, $f$ is uniquely determined by the $10$ numbers $\{f(x_i)\}_{i=1}^4,\{(\Delta+\lambda)f(y_i)\}_{i=1}^6$, which contradicts the fact that $l(\Gamma)$ is $11$ dimensional. \hfill$\square$\\

Now, we return to investigate the Sierpinski lattice $\widetilde{\SG}$, applying the above two lemmas locally. 
\begin{lemma}\label{Q1lemma3}
	Let $\lambda\notin\{1,2,4,5,6\}$, and consider $\Delta_{(-k)}:\ell^p(\lSG^{(-k)})\to  \ell^p(\lSG^{(-k)})$, $1\leq p\leq \infty$. Then $\lambda+\Delta_{(-k)}$ is invertible if and only if $R(\lambda)+\Delta_{(-k-1)}$ is invertible.
\end{lemma}

\noindent\textit{Proof.} Without loss of generality, we consider the $k=0$ case. Fix any point $y_0\in \lSG^{(-1)}$, and choose a neighbourhood of $y_0$ in $\widetilde{\SG}$ that is isomorphic to $\Gamma$. For each $g\in l(\lSG)$, we define 
\[Tg(y_0)=\sum_{i=0}^6 c_{i,\lambda}g(y_i),\]
where $c_{i,\lambda}$ is defined in Lemma \ref{Q1lemma2}. It is clear that $T$ is bounded from $\ell^p(\lSG)$ to $\ell^p(\lSG^{(-1)})$. In addition, since $c_{0,\lambda}\neq 0$, $T$ is surjective.  

Then by using Lemma \ref{Q1lemma2} at each point of $\lSG^{(-1)}$, the following two systems of equations give the same solutions on $\lSG$,
\begin{equation}\label{Q1eqn2}
(\lambda+\Delta) f(x)=g(x)\text{, on }\lSG,
\end{equation}
and
\begin{equation}\label{Q1eqn3}
\begin{cases}
(\lambda+\Delta) f(x)=g(x),&\text{ for }x\in \lSG\setminus \lSG^{(-1)},\\
(R(\lambda)+\Delta_{(-1)}) f(x)=Tg(x), &\text{ for }x\in \lSG^{(-1)},
\end{cases}
\end{equation}
where $g\in \ell^p(\lSG)$. But (\ref{Q1eqn3}) has a unique solution in $\ell^p(\lSG)$ if and only if 
\begin{equation}\label{Q1eqn4}
(R(\lambda)+\Delta_{(-1)}) f(x)=Tg(x) \text{, on }\lSG^{(-1)},
\end{equation}
has a unique solution in $\ell^{p}(\lSG^{(-1)})$. The lemma follows immediately from the equivalence of solvabilty  and uniqueness of solutions to (\ref{Q1eqn2}) and (\ref{Q1eqn4}).\hfill$\square$\\

In fact, in Lemma \ref{Q1lemma3}, the only exceptions are $\lambda=\{2,6\}$. The cases for $\lambda=1,4,5,6$ are easy to check, while the case $\lambda=2$ needs a little more work. Luckily, by a similar idea as the proof of Lemma \ref{Q1lemma1}, \ref{Q1lemma2} and \ref{Q1lemma3}, we will get the following lemma. 

\begin{lemma}\label{Q1lemma4}
	Let $\lambda\notin\{1,3,4,5,6\}\cup R^{\circ-1}\{1,2,5\}$, and consider $\Delta_{(-k)}:\ell^p(\lSG^{(-k)})\to  \ell^p(\lSG^{(-k)})$,$1\leq p\leq \infty$. Then $(\lambda+\Delta_{(-k)})^{-1}$ is invertible if and only if $(R^{\circ 2}(\lambda)+\Delta_{(-k-2)})^{-1}$ is invertible.
	
	In particular, $(2+\Delta_{(-k)})^{-1}$ is invertible if and only if $(-6+\Delta_{(-k-2)})^{-1}$ is invertible.
\end{lemma} 

We end this section with the proof of Theorem \ref{Q1}.\vspace{0.1cm}

\noindent\textit{Proof of Theorem \ref{Q1}.} Clearly, $\Sigma_6=\{6\}\cup  (\bigcup_{m=0}^\infty R^{\circ -m}\{3\})\subset \sigma(\Delta)$. See \cite{T} for the eigenfunctions for $\lambda\in \Sigma_6$. Thus, it is easy to see that 
\[ \mathcal{J}\cup \Sigma_6=\overline{\Sigma_6}\subset \sigma(\Delta).\] 

It remains to show that $\sigma(\Delta)\subset\mathcal{J}\cup \Sigma_6$. It suffices to show that if $\lambda \notin \mathcal{J}\cup \Sigma_6$, then $\lambda+\Delta$ is invertible. We consider two cases below.

First, consider $\lambda\notin \mathcal{J}\cup \Sigma_6\cup (\bigcup_{m=0}^\infty R^{\circ-m}\{2\})$. It is easy to see from (\ref{eqn11}) that 
\[\|\Delta_{(-k)}\|_{op}\leq 8,\]
where $\|\cdot\|_{op}$ stands for the operator norm of $\Delta_{(-k)}:\ell^p(\lSG^{(-k)})\to  \ell^p(\lSG^{(-k)})$. On the other hand, by the definition of $\mathcal{J}$, there exists $k\geq 0$ such that 
\[|R^{\circ k}(\lambda)|>8\geq \|\Delta_{(-k)}\|_{op},\]
which implies that $R^{\circ k}(\lambda)+\Delta_{(-k)}$ is invertible. By using Lemma \ref{Q1lemma3} repeatedly, we see that $\lambda+\Delta$ is invertible.

Second, if $\lambda\in \bigcup_{m=0}^\infty R^{\circ-m}\{2\}$, by using Lemma \ref{Q1lemma3} and \ref{Q1lemma4}, and a same argument as the first case, we can show that $\lambda+\Delta$ is also invertible.\hfill$\square$  

\section{A decomposition of the $\ell^p$-spectrum on $\lSG$ }

In this section, we focus on characterizing each point in the $\ell^p$-spectrum. We will point out that the $1<p<\infty$ case and the $p=1$ case are very different. A full description of the spectral points in the $1<p<\infty$ case is easy with the method developed by A. Teplyaev \cite{T}, while the $p=1$ case is much complicated and we can only give a full answer for the lattices with a boundary point. 

As preparation, we define the \textit{inner product} of real functions on $\lSG$ as follows, $$<f,g>=\sum_{x\in \lSG}\mu_xf(x)g(x),$$
where
$\mu_x=\begin{cases}
1, &\text{ if  $x$  is not a boundary point}, \\
1/2, &\text{ if  $x$  is a boundary point}.
\end{cases}$

\begin{lemma}\label{Q2lemma1}
	Let $f\in \ell^p(\lSG)$ and $g\in \ell^q(\lSG)$ with $\frac{1}{p}+\frac{1}{q}=1$, $1\leq p<\infty$, then we have $<\Delta f,g>=<f,\Delta g>$.
\end{lemma}
\noindent\textit{Proof.} It is easy to see that the summation below converges absolutely, so we can rearrange the order, 
\[\sum_{x\in \lSG}\sum_{y:y\sim x}g(x)f(y)=\sum_{x\in \lSG}\sum_{y:y\sim x}f(x)g(y).\]
As a result,
\[\quad\qquad\qquad\begin{aligned}
<\Delta f,g>&=\sum_{x\in \lSG} \mu_x g(x)\sum_{y:y\sim x}\frac{1}{\mu_x}\big(f(y)-f(x)\big)\\
&=-4<f,g>+\sum_{x\in \lSG}\sum_{y:y\sim x}g(x)f(y)\\
&=-4<f,g>+\sum_{x\in \lSG}\sum_{y:y\sim x}f(x)g(y)=<f,\Delta g>.\qquad\qquad\qquad\square
\end{aligned}\]

In the following, we use $\sigma_c(\Delta)$ to denote the \textit{continuous spectrum} of the Laplacian, and $\sigma_p(\Delta)$ for the \textit{point spectrum}, $\sigma_r(\Delta)$ for the \textit{residue spectrum}. More precisely,
\[
\begin{cases}
\sigma_p(\Delta)=\{\lambda\in \mathbb{C}:(\lambda+\Delta) \text{ is not injective}\},\\
\sigma_c(\Delta)=\{\lambda\in \mathbb{C}\setminus \sigma_p(\Delta):(\lambda+\Delta) \text{ is not surjective, }(\lambda+\Delta)(\ell^p(\lSG))\text{ is dense in }\ell^p(\lSG)\},\\
\sigma_r(\Delta)=\{\lambda\in \mathbb{C}\setminus \sigma_p(\Delta):(\lambda+\Delta) \text{ is not surjective, }(\lambda+\Delta)(\ell^p(\lSG))\text{ is not dense in }\ell^p(\lSG)\},
\end{cases}
\]
noticing that $\sigma_p(\Delta)$, $\sigma_c(\Delta)$ and $\sigma_r(\Delta)$ actually depend on $1\leq p<\infty$, although we do not include $p$ in the notations. Readers can refer to the book \cite{RS} (Chapter 6) or \cite{Rudin} (page 343) for a discussion on the different spectrums. In particular, the residue spectrum is always empty for self-adjoint operators on Hilbert spaces. 

For self-adjoint operators on Hilbert spaces, a stronger description of the spectrum can be given with the spectral measure. In particular, one can divide the measure into three parts (see books \cite{RS,Rudin}): absolutely continuous spectrum, singularly continuous spectrum, and pure point spectrum. In \cite{T}, Teplyaev showed the spectrum of the Neumann $\ell^2$-Laplacian on $\lSG$ is always pure point. The existence of singularly continuous spectrum of Laplacian on $\lSG$ and on some other fractal lattices was also proved, see \cite{CT,Q}. We have to point out that this decomposition applies to the spectral measure, rather than a point to point classification we are considering in this paper. 

As a consequence of Lemma \ref{Q2lemma1}, we have the following criterion for $\lambda$ to be a residue spectral point.

\begin{lemma}\label{Q2lemma2}
Let $1\leq p<\infty$, $\frac{1}{p}+\frac{1}{q}=1$, and assume $\lambda\notin\sigma_p(\Delta)$ for $\Delta:\ell^p(\lSG)\to \ell^p(\lSG)$.  Then $\lambda\in \sigma_r(\Delta)$ for $\Delta:\ell^p(\lSG)\to \ell^p(\lSG)$ if and only if $\lambda\in \sigma_p(\Delta)$ for $\Delta:\ell^q(\lSG)\to \ell^q(\lSG)$.
\end{lemma}

\noindent\textit{Proof.} The lemma is an easy application of Lemma \ref{Q2lemma1}. In fact, if $(\lambda+\Delta)(\ell^p(\lSG))$ is not dense in $\ell^p(\lSG)$, there exists a non-zero $f_\lambda\in \ell^q(\lSG)$ such that $$<f,(\lambda+\Delta)f_\lambda>=<(\lambda+\Delta)f,f_\lambda>=0$$ for any $f\in \ell^p(\lSG)$. This shows that $(\lambda+\Delta)f_\lambda=0$. Conversely, it is clear that if $\lambda$ is an eigenvalue of $\Delta:\ell^q(\lSG)\to \ell^q(\lSG)$ with a corresponding eigenfunction $f_\lambda\in \ell^q(\lSG)$, then $<(\lambda+\Delta)f,f_\lambda>=0$ for any $f\in \ell^p(\lSG)$. \hfill$\square$  \vspace{0.1cm}

Lemma 3.2 is an easy consequence of the fact that $\Delta:\ell^q(\lSG)\to \ell^q(\lSG)$ is the conjugacy operator of $\Delta:\ell^p(\lSG)\to \ell^p(\lSG)$. Readers can find a general result about operators on Banach spaces in \cite{RS} (the proposition above Theorem \uppercase\expandafter{\romannumeral6}.8).

\subsection{Lattices with a boundary point}
In this part, we will characterize each spectral point for $\Delta:\ell^p(\lSG)\to\ell^p(\lSG)$, given the condition that $\lSG$ is a lattice with a boundary point. The result is stated as follows.

\begin{theorem}\label{Q2thm3}
Write 
$$\Sigma_4=\bigcup_{m=0}^\infty R^{\circ-m}\{4\},\quad \Sigma_5=\bigcup_{m=0}^\infty R^{\circ-m}\{5\} \quad \text{and} \quad  \Sigma_6=\{6\}\bigcup\big(\bigcup_{m=0}^\infty R^{\circ-m}\{3\}\big).$$
	
(a). For $\Delta:\ell^p(\lSG)\to \ell^p(\lSG)$ with $1<p<\infty$, we have $\sigma_p(\Delta)=\Sigma_5\cup\Sigma_6$ and $\sigma_c(\Delta)=\mathcal{J}\setminus \Sigma_5$. There is no residue spectral point.

(b). For $\Delta:\ell^1(\lSG)\to \ell^1(\lSG)$, we have all three types of spectral points as follows,
\[\sigma_p(\Delta)=\Sigma_5\cup\Sigma_6,\quad\sigma_c(\Delta)=\mathcal{J}\setminus (\Sigma_4\cup \Sigma_5\cup\{0\}),\quad \sigma_r(\Delta)=\{0\}\cup\Sigma_4.\]
\end{theorem}

In the following Proposition \ref{Q2prop4}, we will see that for the lattice $\lSG$ with a boundary point, any eigenfunction $f$ to $\lambda\notin \{0\}\cup\Sigma_4\cup\Sigma_5\cup\Sigma_6$ is unbounded. Theorem \ref{Q2thm3} is then an immediate consequence by Lemma \ref{Q2lemma2}.

\begin{figure}[htp]
\includegraphics[width=5cm]{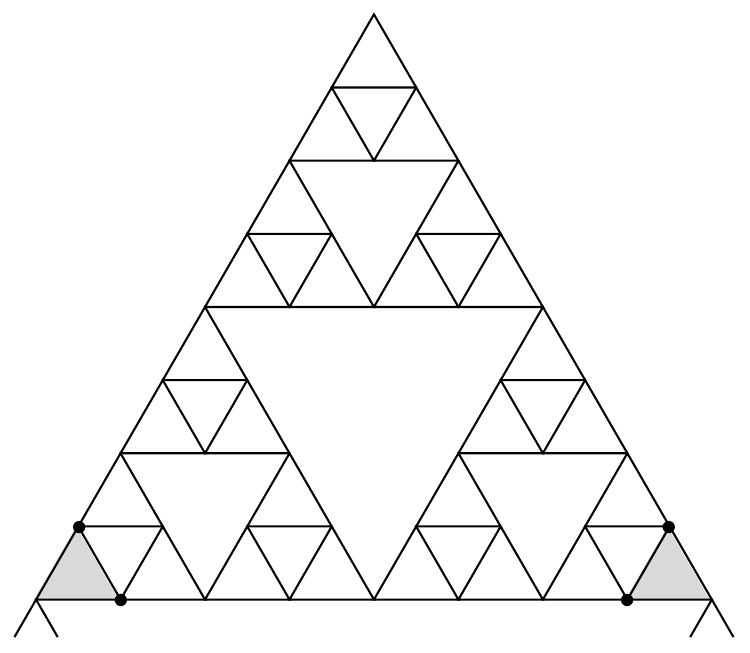}
\caption{The set of vertices $\{F_0^{-m}F_i^{m}q_j:i,j=1,2\}$ with $m=3$.}\label{explode}
\end{figure}

In particular, we will apply the fact that $\lSG$ is expanding towards a fixed direction  if there exists a boundary point, say $F_0^{-1}$ without loss of generality. By some computations (using the decimation algorithm, Proposition \ref{decimation} and more detail in \cite{S1}, Lemma 3.2.1), we will see that $f$ is unbounded on the vertices $\bigcup_{m\geq 1}\{F_0^{-m}F_i^{m}q_j:i,j=1,2\}$, which are  neighbours of the boundaries of $F_1^{-m}V_m$. See Figure \ref{explode} for an illustration, where $\{F_0^{-3}F_i^{3}q_j:i,j=1,2\}$ are dotted, and the triangles containing the points are shaded.

\begin{proposition}\label{Q2prop4}
  For $\Delta:\ell^\infty(\lSG)\to \ell^\infty(\lSG)$, we have $\sigma_p(\Delta)=\{0\}\cup\Sigma_4\cup\Sigma_5\cup\Sigma_6$. 
\end{proposition} 

\noindent\textit{Proof.} Without loss of generality, we take $\omega=0000\cdots$ and $\lSG=\bigcup_{m=0}^\infty F_0^{-m}V_m$, since any two lattices with a boundary point are isomorphic to each other \cite{T}. The following proof relies on Lemma \ref{Q2lemma5} and Lemma \ref{Q2lemma6}, which will be stated later.

The existence of $5$-series and $6$-series eigenfunctions is a well-known result of the spectral decimation, see \cite{T}. See Appendix for the existence of the $4$-series eigenfunctions. So it suffices to show that there are no other $\ell^\infty$-eigenvalues.

Take $\lambda\in \mathcal{J}\setminus \big(\{0\}\cup\Sigma_4\cup\Sigma_5\big)$, and let $f$ be a $\lambda$-eigenfunction. For convenience, we write $q_i^{(-m)}=F_0^{-m}q_i,i=0,1,2$. By direct computation, and using the decimation method, we get
\begin{equation}\label{eqn31}
\begin{aligned}
\begin{pmatrix}f(q_0)\\f(q_1)\\f(q_2)\end{pmatrix}&=
\begin{pmatrix}
1 & 0 & 0\\
\frac{4-\lambda}{(2-\lambda)(5-\lambda)} & \frac{4-\lambda}{(2-\lambda)(5-\lambda)} & \frac{2}{(2-\lambda)(5-\lambda)}\\
\frac{4-\lambda}{(2-\lambda)(5-\lambda)} & \frac{2}{(2-\lambda)(5-\lambda)} & \frac{4-\lambda}{(2-\lambda)(5-\lambda)}
\end{pmatrix}
\begin{pmatrix}f(q_0)\\f(q^{(-1)}_1)\\f(q^{(-1)}_2)\end{pmatrix}
\\
&=f(q_0)\begin{pmatrix}1\\1-\frac{\lambda}{4}\\1-\frac{\lambda}{4}\end{pmatrix}
+\frac{\lambda\big(f(q^{(-1)}_1)-f(q^{(-1)}_2)\big)}{2R(\lambda)}
\begin{pmatrix}0\\1\\-1\end{pmatrix}\\
&+\big(\frac{f(q_1^{(-1)})+f(q_2^{(-1)})}{2}-(1-\frac{R(\lambda)}{4})f(q_0)\big)\frac{6-\lambda}{(2-\lambda)(5-\lambda)}
\begin{pmatrix}0\\1\\1\end{pmatrix},
\end{aligned}
\end{equation}
as $f$ is an $R(\lambda)$-eigenfunction on $\lSG^{(-1)}$. In addition, 

\begin{equation}\label{eqn32}
f(q_1^{(-m)})+f(q_2^{(-m)})-(2-\frac{R^{\circ m}(\lambda)}{2})f(q_0)=0, \forall m\geq 0.
\end{equation}
By using (\ref{eqn31}) and (\ref{eqn32}) repeatedly and using symmetry, we get 
\begin{equation}\label{eqn33}
\begin{pmatrix}f(q_0)\\f(q_1)\\f(q_2)\end{pmatrix}=
f(q_0)\begin{pmatrix}1\\1-\frac{\lambda}{4}\\1-\frac{\lambda}{4}\end{pmatrix}
+\frac{\lambda\big(f(q^{(-m)}_1)-f(q^{(-m)}_2)\big)}{2R^{\circ m}(\lambda)}
\begin{pmatrix}0\\1\\-1\end{pmatrix},
\end{equation} 
and 
\begin{equation}\label{eqn34}
\begin{aligned}
\begin{pmatrix}f(F_0^{-m}F^m_1q_0)\\f(q^{(-m)}_1)\\f(F_0^{-m}F^m_1q_2)\end{pmatrix}=
f(q^{(-m)}_1)\begin{pmatrix}1-\frac{\lambda}{4}\\1\\1-\frac{\lambda}{4}\end{pmatrix}
+\frac{\lambda\big(f(q_0)-f(q^{(-m)}_2)\big)}{2R^{\circ m}(\lambda)} 
\begin{pmatrix}1\\0\\-1\end{pmatrix}\\
+\big(\frac{f(q_0)+f(q_2^{(-m)})}{2}-(1-\frac{R^{\circ m}(\lambda)}{4})f(q_1^{(-m)})\big)\prod_{l=0}^{m-1}\frac{6-R^{\circ l}(\lambda)}{(2-R^{\circ l}(\lambda))(5-R^{\circ l}(\lambda))}\begin{pmatrix}1\\0\\1\end{pmatrix}.
\end{aligned}
\end{equation} 
As a consequence of (\ref{eqn32}) and (\ref{eqn33}), we get the estimate
\begin{equation}\label{eqn35}
\max\{|f(q_0)|,|f(q_1)|,|f(q_2)|\}\leq C (1+\frac{\lambda}{R^{\circ m}(\lambda)})\cdot\max\{|f(q^{(-m)}_1)|,|f(q^{(-m)}_2)|\}. 
\end{equation}
On the other hand, by using (\ref{eqn32}) and (\ref{eqn34}), we get the equation
\begin{equation}
\begin{aligned}
f(F_0^{-m}F^m_1q_0)+f(F_0^{-m}&F^m_1q_2)=(2-\frac{\lambda}{2})f(q_1^{(-m)})\\
+\big((\frac{2}{4-R^{\circ m}(\lambda)}&-2+\frac{R^{\circ m}(\lambda)}{2})f(q_1^{(-m)})+(\frac{2}{4-R^{\circ m}(\lambda)}+1)f(q_2^{(-m)})\big)P_m\\
&=(2-\frac{\lambda}{2})f(q_1^{(-m)})+P_m(\lambda)(a_mf(q_1^{(-m)})+b_mf(q_2^{(-m)})).
\end{aligned}
\end{equation}
where 
\begin{equation}\label{pm}
P_m(\lambda)=\prod_{l=0}^{m-1}\frac{6-R^{\circ l}(\lambda)}{(2-R^{\circ l}(\lambda))(5-R^{\circ l}(\lambda)},
\end{equation}
 $a_m=\frac{2}{4-R^{\circ m}(\lambda)}-2+\frac{R^{\circ m}(\lambda)}{2}$ and $b_m=\frac{2}{4-R^{\circ m}(\lambda)}+1$.
By symmetry, we also have 
\begin{equation}
\begin{aligned}
f(F_0^{-m}F^m_2q_0)&+f(F_0^{-m}F^m_2q_1)\\
&=(2-\frac{\lambda}{2})f(q_2^{(-m)})+P_m(\lambda)(a_mf(q_2^{(-m)})+b_mf(q_1^{(-m)})).
\end{aligned}
\end{equation}

It is easy to check that $a_m^2\neq b_m^2$ if $R^{\circ m}(\lambda)\in \mathcal{J}\setminus \{0,4\}$. By using Lemma \ref{Q2lemma5} below, we can find an increasing sequence $\{m_k\}_{k=1}^\infty$ such that $R^{\circ m_k}(\lambda)$ is bounded away from $0$ and $4$. Thus, we have 
\begin{equation}\label{eqn38}
\begin{aligned}
&\max\{|f(F_0^{-m_k}F^{m_k}_1q_0)+f(F_0^{-m_k}F^{m_k}_1q_2)|,|f(F_0^{-m_k}F^{m_k}_2q_0)+f(F_0^{-m_k}F^{m_k}_2q_1)|\}\\
\geq & (2C |P_{m_k}(\lambda)|-|2-\frac{\lambda}{2}|)\max \{|f(q^{(-m_k)}_1)|,|f(q^{(-m_k)}_2)|\},
\end{aligned}
\end{equation}
where $C$ is independent of $k$. However, according to Lemma \ref{Q2lemma6}, we can see that $\lim\limits_{k\to\infty}|P_{m_k}(\lambda)|=+\infty$. Combining the estimates (\ref{eqn35}) and (\ref{eqn38}), and letting $k\to\infty$, we see that $f$ is unbounded. \hfill$\square$\vspace{0.2cm}

At the end of this subsection, we prove the lemmas that are used in the proof of Proposition \ref{Q2prop4}.

First, let's recall some basic facts about the totally disconnected Julia set, which can be found in many textbooks, see for example \cite{M}. Denote $\varphi_-(x)=\frac{5-\sqrt{25-4x}}{2}$ and $\varphi_+(x)=\frac{5+\sqrt{25-4x}}{2}$.  There is a natural homeomorphism $\pi$ from the Cantor set $\mathcal{C}=\{-,+\}^{\mathbb{N}}$ to $\mathcal{J}$ defined as follows 
\[\pi(\eta)=\bigcap_{m=1}^\infty \varphi_{\eta_1}\varphi_{\eta_2}\cdots \varphi_{\eta_m}(\mathcal{J}),\]
for $\eta=\eta_1\eta_2\cdots\in \mathcal{C}$. In particular, $\pi(---\cdots )=0$ and $\pi(+++\cdots)=4$. In addition, define the left shift operator  $\iota: \mathcal{C}\rightarrow \mathcal{C}$ by $\iota(\eta_1\eta_2\eta_3\cdots)=\eta_2\eta_3\cdots$. Then
\[R\circ \pi(\eta)=\pi\circ \iota (\eta), \forall \eta\in \mathcal{C}.\]

\begin{lemma}\label{Q2lemma5}
	Let $\lambda\in \mathcal{J}\setminus \big(\{0\}\cup\Sigma_4\cup\Sigma_5\big)$.
	
	(a). There is a sequence $m_1<m_2<\cdots $ such that $\min\{|R^{\circ m_k}(\lambda)|,|R^{\circ m_k}(\lambda)-4|\}>C$, where $C$ is a positive constant independent of $\lambda$. 
	
	(b). There is an infinite sequence $n_1<n_2<\cdots $ such that $R^{\circ n_k}(\lambda)\in \varphi_-(\mathcal{J})$.
\end{lemma}
\noindent\textit{Proof.} (a). Clearly, we can find countably infinite different positive integers $m$ such that $\iota^{\circ m}\pi^{-1}(\lambda)$ are of the form $+-\cdots$ or $-+\cdots$. This means $R^{\circ m}(\lambda)\in \varphi_{-}\varphi_{+}(\mathcal{J})\cup \varphi_{+}\varphi_{-}(\mathcal{J})$, and thus $R^{\circ m}(\lambda)$ is bounded away from $\{0,4\}$.

(b). Clearly, we can find countably infinite different integers $n$ such that $\big(\iota^{\circ n}\pi^{-1}(\lambda)\big)_1=-$.\hfill$\square$\vspace{0.2cm}

Next, we give an estimate for $P_m(\lambda)$ in $(\ref{pm})$.

\begin{lemma}\label{Q2lemma6}  For $\lambda\in \mathcal{J}\setminus \Sigma_4$, we have 
	$$\lim_{m\to\infty} |P_m(\lambda)|=\lim\limits_{m\to\infty}\prod_{l=0}^{m-1}|\frac{6-R^{\circ l}(\lambda)}{(2-R^{\circ l}(\lambda))(5-R^{\circ l}(\lambda)}|=\infty.$$
\end{lemma}

\noindent\textit{Proof.} Using the fact that $R(x)=x(5-x)$ and by direct computation, we can show
\[|P_{m}(\lambda)|=\big|\frac{\lambda(6-\lambda)}{R^{\circ m}(\lambda)\big(2-R^{\circ m-1}(\lambda)\big)}\big|\cdot \prod_{l=0}^{m-2}|3-R^{\circ l}(\lambda)|.\]
Thus it suffices to show $\lim\limits_{m\to\infty}\prod_{l=0}^{m}|3-R^{\circ l}(\lambda)|=\infty$.

\begin{figure}[htp]
	\includegraphics[width=12cm]{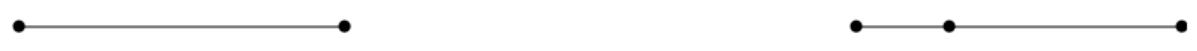}
	\begin{picture}(0,0)
	\put(-300,10){$A$}
	\put(-93,10){$B$}
	\put(-41,10){$C$}
	
	\put(-343,-3){$0$}
	\put(-253,-3){$1.38$}
	\put(-113,-3){$3.62$}
	\put(-77,-3){$4$}
	\put(-10,-3){$5$}
	\end{picture}
	\caption{An illustration of the $A,B,C$ areas and approximated values of the end points.}\label{fig5}
\end{figure}
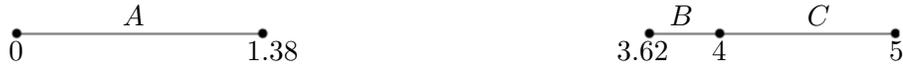

Let $A=\varphi_{-}(\mathcal{J})$, $B=\varphi_{+}(\mathcal{J})\cap (0,4]$ and $C=\varphi_+(\mathcal{J})\cap (4,5]$, so that $\mathcal{J}=A\cup B\cup C$. See Figure \ref{fig5} for an illustration. For $x \in A$, we have $|3-x|>1.5$; for $x \in C$, we have $|3-x|>1$; For $x\in B$, we have $R(x)\in C$ and $(3-x)\big(3-R(x)\big)>1$ by an easy estimate. 

As a consequence, we have the estimates 
\[\begin{aligned}
\prod_{l=0}^{m}|3-R^{\circ l}(\lambda)|&\geq c\big(\prod_{l\in I_{m,A}}|3-R^{\circ l}(\lambda)|\big)\cdot \big(\prod_{l\in I_{m,B}}|(3-R^{\circ l}(\lambda)\big)\big(3-R^{\circ l+1}(\lambda)\big)|\big)\\
&\geq c (\frac{3}{2})^{\# I_{m,A}},
\end{aligned}\]
where $c=\min\{|x-3|:x\in\mathcal{J}\}$, $I_{m,A}=\{0\leq l\leq m:R^{\circ l}(\lambda)\in A\}$ and $I_{m,B}=\{0\leq l\leq m-1:R^{\circ l}(\lambda)\in B\}$. The lemma follows immediately from Lemma \ref{Q2lemma5} (b).\hfill$\square$

\subsection{Lattices with no boundary}
In A. Teplyaev's work \cite{T}, it was shown that the localized eigenfunctions, i.e. $5$-series or $6$-series functions, form a complete basis of the $\ell^2(\lSG)$ space. The basic idea is to find a localized eigenfunction $f_\lambda$ such that $<f,f_\lambda>\neq 0$ for each nonzero $f\in \ell^2(\lSG)$. The proof can be easily extended to $C_0(\lSG)$ case, where $$C_0(\lSG)=\{f\in \ell^\infty(\lSG):\lim_{x\to\infty} f(x)=0\}.$$

\begin{lemma}\label{Q2lemma7}
	For any nonzero $f\in C_0(\lSG)$, there exists a localized eigenfunction $f_\lambda$ of $\Delta$ such that $<f_\lambda,f>\neq 0$.
\end{lemma}
\begin{proof}
Recall that $\lSG=\bigcup_{m=1}^\infty F_{\omega_1}^{-1}F_{\omega_2}^{-1}\cdots F_{\omega_m}^{-1}(V_m)$ for some $\omega\in \{0,1,2\}^\infty$. In addition, since $\lSG$ has no boundary point, there are infinitely many $m\in \mathbb{N}$ such that $\omega_m\neq \omega_{m+1}$ by Lemma 2.3 of \cite{T}. 

For short, we write $V_{-m}=F_{\omega_1}^{-1}F_{\omega_2}^{-1}\cdots F_{\omega_m}^{-1}(V_m)$ for $m\geq 1$. Let $\mathcal S_3$ be the permutation group on $\{0,1,2\}$. For each $\sigma\in \mathcal S_3$, we define the local symmetry $T_{m,\sigma}:V_{-m}\to V_{-m}$ by
\[T_{m,\sigma}(F_{\omega_1}^{-1}F_{\omega_2}^{-1}\cdots F_{\omega_m}^{-1}F_{w_1}F_{w_2}\cdots F_{w_m}q_i)=F_{\omega_1}^{-1}F_{\omega_2}^{-1}\cdots F_{\omega_m}^{-1}F_{\sigma(w_1)}F_{\sigma(w_2)}\cdots F_{\sigma(w_m)}q_{\sigma(i)},\]
for any $w=w_1w_2\cdots w_m\in \{0,1,2\}^m$. The space $S_m=\{f\in l(\lSG):f(x)=0,\forall x\notin V_{-m}, \text{ and }f\circ T_{m,\sigma}=sgn(\sigma)f, \forall \sigma\in \mathcal S_3\}$, where $sgn(\sigma)$ is the sign of the permutation $\sigma$, is spanned by localized eigenfunctions.
 
Let $x\in \lSG$ such that $f(x)\neq 0$, and choose $m\geq 1$ so that $\omega_m\neq \omega_{m+1}$, $x\in V_{-m+1}$, and $|f(y)|\leq\frac{1}{6}|f(x)|$, $\forall y\notin V_{-m+1}$. Define $g\in S_{m+1}$ as
\[g=\sum_{\sigma\in \mathcal S_3}sgn(\sigma)1_{T_{m+1,\sigma}(x)},\]
where $1_{z}$ is the indicator function, i.e. $1_z(z)=1$ and $1_z(y)=0$ if $y\neq z$. Clearly, if $\sigma$ is not the identity permutation, we have $T_{m+1,\sigma}(x)\notin V_{-m+1}$, so $|<g,f>|=|\sum_{\sigma\in \mathcal S_3}sgn(\sigma)f(T_{m+1,\sigma}(x))|\geq |f(x)|-\frac{5}{6}|f(x)|>0$. The lemma follows immediately as $g$ is a linear combination of localized eigenfunctions. 
\end{proof}

As an immediate consequence, there is no eigenvalue of $\Delta$ on $C_0(\lSG)$ other than the $5$ or $6$ series.

\begin{proposition}\label{Q2prop8}
	Let $\Delta:C_0(\lSG)\to C_0(\lSG)$, we have $\sigma_p(\Delta)=\Sigma_5\cup\Sigma_6$.
\end{proposition}
\noindent\textit{Proof.} Assume there exists an eigenvalue $\lambda\notin \Sigma_5\cup\Sigma_6$, and let $f_{\lambda}$ be a corresponding eigenfunction. By Lemma \ref{Q2lemma7}, there is a localized eigenfunction $f_{\lambda'}$ such that $<f_\lambda,f_{\lambda'}>\neq 0$. Noticing that  $\lambda\neq 0$ and $\lambda'\in\Sigma_5\cup\Sigma_6$, we have
\[\lambda^{-1}<-\Delta f_\lambda,f_{\lambda'}>=<f_\lambda, f_{\lambda'}>={\lambda'}^{-1}< f_\lambda,-\Delta f_{\lambda'}>.\]
Then by Lemma \ref{Q2lemma1},  this implies that $<f_\lambda,f_{\lambda'}>=0$,  a contradiction. \hfill$\square$
\begin{theorem}
	 For $1<p<\infty$, $\Delta:\ell^p(\lSG)\to \ell^p(\lSG)$ has point spectrum $\sigma_p(\Delta)=\Sigma_5\cup\Sigma_6$, and continuous spectrum $\sigma_c(\Delta)=\mathcal{J}\setminus \Sigma_5$. There is no residue spectrum.
\end{theorem}
\noindent\textit{Proof.} As a direct consequence of Proposition \ref{Q2prop8}, $\sigma_p(\Delta)=\Sigma_5\cup\Sigma_6$. In addition, we see that there is no residue spectral point by Lemma \ref{Q2lemma2}.\hfill$\square$ \vspace{0.2cm}

However, the $\ell^1$-spectrum of the Laplacian is much complicated on the lattices without boundary, as the approach for Proposition \ref{Q2prop4} is no longer available in this situation. It seems possible that the eigenvalues for $\Delta:\ell^\infty(\lSG)\to \ell^\infty(\lSG)$ depend on the generating sequence $\omega$ of the lattice. We leave the problem for further study.\vspace{0.2cm}

\noindent\textbf{Problem}. \textit{Describe all the eigenvalues of $\Delta:\ell^\infty(\lSG)\to \ell^\infty(\lSG)$. In particular, whether the set of $\ell^\infty$-eigenvalues depends on the generating sequence $\omega\in \{0,1,2\}^\infty$?}

\subsection{Infinite Sierpinski lattice of Barlow and Perkins}
\begin{figure}[htp]
	\includegraphics[width=7cm]{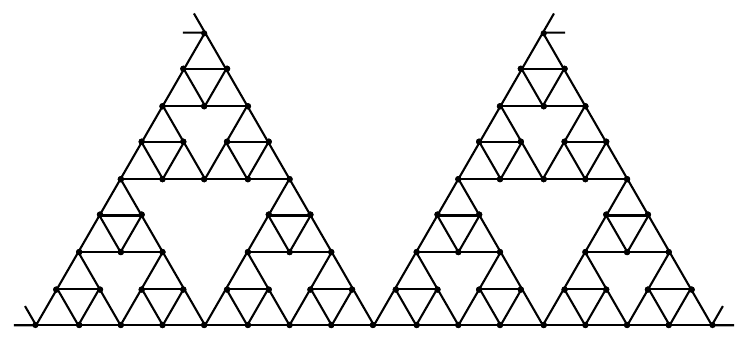}
	\caption{The infinite Sierpinski lattice $\lSG_{BP}$ of Barlow and Perkins.}\label{bplattice}
\end{figure}

At the end of this section, we point out that there is another important infinite Sierpinski lattice in history, introduced by M.T. Barlow and E.A. Perkins \cite{BP}. We write $\lSG_{BP}$ for this lattice, see Figure \ref{bplattice}.

This example is not an infinite Sierpinski lattice considered above, it is the union of two copies of $\lSG$ which are joined at the boundary. By symmetry, one can easily see that the Laplacian on this lattice is equivalent to a direct sum of $\Delta$ and $\Delta_D$, where $\Delta$ is the Neumann Laplacian we consider above, and $\Delta_D$ is the Dirichlet Laplacian on $\lSG$ with a boundary point. 

The argument in Section 2 applies here without difficulty. In addition, one can modify the argument of Proposition \ref{Q2prop4} to show that the $\ell^p$-eigenvalues of $\Delta$ on $\lSG_{BP}$ is still $\Sigma_5\cup\Sigma_6$. In fact, we have to drop (\ref{eqn32}) in the argument since we no longer have the Neumann boundary condition on a half of $\lSG_{BP}$, but still we can see if $\lambda\in \mathcal{J}\setminus \big(\{0\}\cup\Sigma_4\cup\Sigma_5\big)$ and $f$ is a non-zero $\lambda$-eigenfunction, $|f|$ is bounded below on infinitely many points of the form $F_0^{-m}F_i^mq_j$, with $i\in \{1,2\}$ and $j\neq i$. We leave the detail to readers. 

\begin{theorem}\label{bpthm}
	For $1<p<\infty$, $\Delta:\ell^p(\lSG_{BP})\to \ell^p(\lSG_{BP})$ has point spectrum $\sigma_p(\Delta)=\Sigma_5\cup\Sigma_6$, and continuous spectrum $\sigma_c(\Delta)=\mathcal{J}\setminus \Sigma_5$. There is no residue spectrum.
\end{theorem}

On the other hand, the argument in Proposition \ref{Q2prop4}, showing that the eigenfunctions are not bounded, relies heavily on (\ref{eqn32}). It is not clear how to overcome this, and we do not know whether the same result of Theorem \ref{Q2thm3} (b) holds on $\lSG_{BP}$.

\section{Appendix}
In this appendix, we construct the $4$-eigenfunctions on $\widetilde{\SG}$. Note that this induces a class of eigenvalues $\Sigma_4=\bigcup_{m=0}^\infty R^{\circ-m}\{4\}$.

We introduce the following orthogonal matrices
\begin{equation}\label{extend}
A_1=
\begin{pmatrix}
1&0&0\\
0&0&-1\\
0&-1&0
\end{pmatrix},\quad
A_2=
\begin{pmatrix}
0&0&-1\\
0&1&0\\
-1&0&0
\end{pmatrix},
\quad
A_3=
\begin{pmatrix}
0&-1&0\\
-1&0&0\\
0&0&1
\end{pmatrix}.
\end{equation}\quad

Recall that if we fix an infinite word $\omega=\omega_1\omega_2\cdots$, then there is a Sierpinski lattice defined by $\lSG=\bigcup_{m=0}^\infty F_{\omega_1}^{-1}F^{-1}_{\omega_2}\cdots F^{-1}_{\omega_m} V_m$. 
For convenience, we write $$q_i^{(-m)}=F_{\omega_1}^{-1}F^{-1}_{\omega_2}\cdots F^{-1}_{\omega_m}(q_i),\quad i=0,1,2,$$
and  write 
$$q_{li}^{(-m)}=F_{\omega_1}^{-1}F^{-1}_{\omega_2}\cdots F^{-1}_{\omega_m}(F_lq_i),\quad l\in W_m=\{0,1,2\}^m.$$
Clearly $q_i=q_{\omega_m\omega_{m-1}\cdots\omega_1 i}^{(-m)}$. 

\begin{proposition}
	(a). If $\lSG$ has no boundary, then there is a three dimensional $\ell^\infty$-eigenspace of $\Delta$ corresponding to $4$. 
	
	(b). If $\lSG$ has a boundary point, then there is a two dimensional $\ell^\infty$-eigenspace of $\Delta$ corresponding to $4$.
\end{proposition}

\begin{figure}[htp]
	\includegraphics[width=4cm]{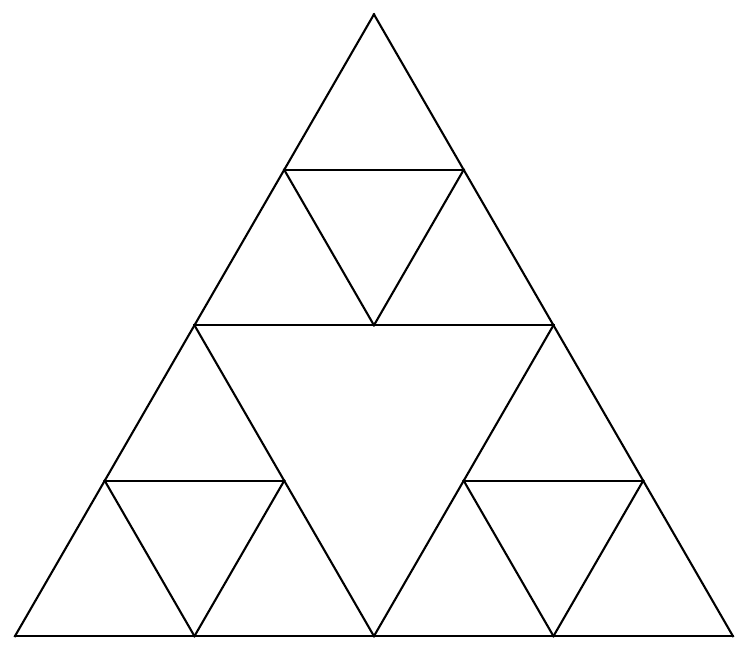}\qquad
	\includegraphics[width=4cm]{gamma2.pdf}
	\begin{picture}(0,0)
	\put(-73,24){$a$}
	\put(-91,-6){$b$}
	\put(-63,-6){$c$}
	
	\put(-130,-1){$-a$}
	\put(-103,48){$-b$}
	\put(-116,24){$-c$}
	
	\put(-4,-1){$b$}
	\put(-67,98){$-c$}
	\put(-32,48){$a$}
	
	\put(-4,-1){$b$}
	\put(-67,98){$-c$}
	\put(-32,48){$a$}
	
	\put(-89,71){$-a$}
	\put(-46,71){$b$}
	\put(-63,43){$c$}
	
	\put(-42,-6){$-a$}
	\put(-61,24){$-b$}
	\put(-22,24){$-c$}
	
	\put(-210,26){$q_0$}
	\put(-228,-5){$q_1$}
	\put(-198,-5){$q_2$}
	\end{picture}
	\caption{An illustration for extending $f$ to be a $4$-eigenfunction.(We take $\omega_1=2,\omega_2=1$ as shown in the left picture.) }\label{fig6}
\end{figure}

\noindent\textit{Proof.} (a). Let $f(q_0)=a,f(q_1)=b,f(q_2)=c$, where $a,b,c$ are arbitrary real number. Define 
\begin{equation}
\begin{pmatrix}
f(q^{(-m)}_0)\\
f(q^{(-m)}_1)\\
f(q^{(-m)}_2)
\end{pmatrix}=
A^{-1}_{\omega_{m}}\cdots A^{-1}_{\omega_{2}}A^{-1}_{\omega_1}\begin{pmatrix}
f(q_0)\\
f(q_1)\\
f(q_2)
\end{pmatrix},
\end{equation}
and 
\begin{equation}\begin{pmatrix}
f(q^{(-m)}_{l0})\\
f(q^{(-m)}_{l1})\\
f(q^{(-m)}_{l2})
\end{pmatrix}=
A_{l_1}A_{l_2}\cdots A_{l_m}\begin{pmatrix}
f(q^{(-m)}_0)\\
f(q^{(-m)}_1)\\
f(q^{(-m)}_2)
\end{pmatrix},\quad\forall l\in W_m.
\end{equation}
See Figure \ref{fig6} for an example of the extension of $f$. One can easily check that $f$ is a $4$-eigenfunction of $\Delta$  on $\widetilde{\SG}$ and $f$ is bounded. By the above construction, we get a three dimensional eigenspace to $4$. 

On the other hand, noticing that $4$ is not a forbidden eigenvalue, a $4$-eigenfunction $f$ is uniquely determined by $f|_{V_0}$. 

(b). The proof of (b) is essentially the same. The eigenspace is $2$ dimensional as a consequence of the eigenvalue equation at the boundary point.\hfill$\square$

\bibliographystyle{amsplain}

\end{document}